\documentclass[a4, 12pt]{amsart}
\usepackage{amssymb}
\usepackage{amstext}
\usepackage{amsmath}
\usepackage{amscd}
\usepackage{amsthm}
\usepackage{latexsym}
\usepackage{amsfonts}
\usepackage{graphicx}
\usepackage{enumerate}

\theoremstyle{plain}
\newtheorem{theorem}{Theorem}[section]

\newtheorem{lemma}[theorem]{Lemma}
\newtheorem{corollary}[theorem]{Corollary}

\theoremstyle{definition}

\theoremstyle{remark}

\newtheorem*{ac}{Acknowledgments}

\numberwithin{equation}{section}

\newtheorem*{tpf}{{\it Proof of Theorem \ref{main}}}
\newtheorem*{spf}{{\it Proof of Corollary \ref{stp}}}

\def\Hom{\operatorname{Hom}}
\def\RHom{\operatorname{{\bf R}Hom}}

\def\Ext{\operatorname{Ext}}
\def\Tor{\operatorname{Tor}}

\def\mod{\operatorname{mod}}

\def\tr{\operatorname{tr}}

\def\tensor{\operatorname{\otimes}}

\def\m{\mathfrak m}

\def\E{\operatorname{E}}

\def\H{\operatorname{H}}

\def\depth{\operatorname{depth}}

\def\pd{\operatorname{pd}}

\def\cpd{\operatorname{\text{$C$}-pd}}

\def\tr{\operatorname{Tr}}

\def\CM{\operatorname{CM}}

\def\add{\operatorname{add}}

\def\Gc{{\mathcal G}}

\def\Ic{{\mathcal I}}

\def\Xc{{\mathcal X}}

\begin{document}

\setlength{\baselineskip}{15pt}

\title[Modules of finite projective dimension]{On the left perpendicular category of the modules of finite projective dimension}

\author{Tokuji Araya}
\address{Liberal Arts Division, Tokuyama College of Technology, Gakuendai, Shunan, Yamaguchi 745-8585, Japan}
\email{araya@tokuyama.ac.jp}

\author{Kei-ichiro Iima}
\address{Department of Liberal Studies, Nara National College of Technology, 22 Yata-cho, Yamatokoriyama, Nara 639-1080, Japan}
\email{iima@libe.nara-k.ac.jp}

\author{Ryo Takahashi}
\address{Department of Mathematical Sciences, Faculty of Science, Shinshu University, 3-1-1 Asahi, Matsumoto, Nagano 390-8621, Japan}
\email{takahasi@math.shinshu-u.ac.jp}

\keywords{perpendicular category, projective dimension, semidualizing module, totally reflexive module, strong test module for projectivity}
\footnote[0]{2010 {\em Mathematics Subject Classification.} 13C60, 13D05, 13H10}
\begin{abstract}
In this paper, we characterize several properties of commutative noetherian local rings in terms of the left perpendicular category of the category of finitely generated modules of finite projective dimension.
As an application we prove that a local ring is regular if (and only if) there exists a strong test module for projectivity having finite projective dimension.
We also obtain corresponding results with respect to a semidualizing module.
\end{abstract}
\maketitle
\section{Introduction}

Throughout this paper, let $R$ be a commutative noetherian local ring with maximal ideal $\m$ and residue field $k$.
All modules considered in this paper are assumed to be finitely generated.

An $R$-module $C$ is said to be {\em semidualizing} if the natural homomorphism $R \to \Hom_R(C,C)$ is an isomorphism and $\Ext_R^i(C,C)=0$ for all $i>0$.
A semidualizing module admits a duality property, which has been defined by Foxby \cite{F}, Vasconcelos \cite{V} and Golod \cite{G}.
A free module of rank one and a dualizing module are semidualizing modules.
The notion of a semidualizing module has been extended by Christensen \cite{C} to a complex which is called a semidualizing complex, and explored in a lot of directions by many authors.
Various homological dimensions with respect to a fixed semidualizing $R$-module $C$ have been invented and investigated (cf. \cite{ATY,C,Ge1,G,TW}).
Among them, the {\em $C$-projective dimension} of a nonzero $R$-module $M$, denoted by $\cpd_RM$, is defined as the infimum of integers $n$ such that there exists an exact sequence of the form
$$
0 \to C^{b_n} \to C^{b_{n-1}} \to \cdots \to C^{b_1} \to C^{b_0} \to M \to 0,
$$
where each $b_i$ is a positive integer.
(The $C$-projective dimension of the zero $R$-module is defined as $-\infty$.)
This is the same thing as $\mathcal{P}_C$-projective dimension; see \cite[Theorems 2.8, 2.11 and Corollary 2.9]{TW}.

An $R$-module $M$ is called {\em totally $C$-reflexive}, where $C$ is a semidualizing $R$-module, if the natural homomorphism $M \rightarrow \Hom_R(\Hom_R(M,C),C)$ is an isomorphism and $\Ext_R^i(M,C)=\Ext_R^i(\Hom_R(M,C),C)=0$ for all $i>0$.
The {\em complete intersection dimension} of $M$, which has been introduced in \cite{AGP}, is defined as the infimum of $\pd_S (M\tensor_RR')-\pd_SR'$ where $R \to R' \leftarrow S$ runs over all quasi-deformations.
Here, a diagram $R \overset{f}{\to} R' \overset{g}{\leftarrow} S$ of homomorphisms of local rings is said to be a quasi-deformation if $f$ is faithfully flat and $g$ is a surjection whose kernel is generated by an $S$-regular sequence.

We denote by $\mod R$ the category of (finitely generated) $R$-modules.
Let $\Gc_{C}(R)$, $\Ic(R)$, and $\add C$ denote the full subcategories of $\mod R$ consisting of all totally $C$-reflexive $R$-modules, consisting of all $R$-modules of complete intersection dimension zero, and consisting of all direct summands of finite direct sums of copies of $C$, respectively.
Let $\Xc_{C}(R)$ be the {\em left perpendicular category} of the category of $R$-modules of finite $C$-projective dimension, that is, the subcategory of $\mod R$ consisting of all $R$-modules $X$ satisfying $\Ext^1_R(X,M)=0$ for each $R$-module $M$ of finite $C$-projective dimension.
We write $\Gc(R)=\Gc_R(R)$ and $\Xc(R)=\Xc_R(R)$.
There are inclusion relations of subcategories of $\mod R$ (cf. Lemma \ref{prop2}):
\begin{gather*}
\Xc(R)\supset\Gc(R)\supset\Ic(R)\supset\add R,\\
\Xc_{C}(R)\supset\Gc_{C}(R)\supset\add C,\ \add R.
\end{gather*}
The main purpose of this paper is to find out what property is characterized by the equalities of $\Xc(R)$ (respectively, $\Xc_{C}(R)$) and each of $\Gc(R)$, $\Ic(R)$, $\add R$ (respectively, each of $\Gc_{C}(R)$, $\add C$, $\add R$).
The main result of this paper is the following theorem.

\begin{theorem}\label{main}
Let $R$ be a commutative noetherian local ring.
\begin{enumerate}[\rm (1)]
\item
The following are equivalent for a semidualizing $R$-module $C$.
\begin{enumerate}[\rm (a)]
\item
$C$ is dualizing.
\item
$\Xc_C(R)=\Gc_C(R)$ holds.
\end{enumerate}
If this is the case, then $R$ is Cohen-Macaulay.
\item
The following are equivalent.
\begin{enumerate}[\rm (a)]
\item
$R$ is Gorenstein.
\item
$\Xc(R)=\Gc(R)$ holds.
\end{enumerate}
\item
The following are equivalent.
\begin{enumerate}[\rm (a)]
\item
$R$ is a complete intersection.
\item
$\Xc(R)=\Ic(R)$ holds.
\end{enumerate}
\item
The following are equivalent.
\begin{enumerate}[\rm (a)]
\item
$R$ is regular.
\item
$\Xc(R)=\add R$ holds.
\item
$\Xc_C(R)=\add C$ holds for some semidualizing $R$-module $C$.
\item
$\Xc_C(R)=\add R$ holds for some semidualizing $R$-module $C$.
\end{enumerate}
\end{enumerate}
\end{theorem}

On the other hand, the notion of a strong test module for projectivity has been introduced and studied by Ramras \cite{R}.
An $R$-module $M$ is called a {\em strong test module for projectivity} if every $R$-module $N$ with $\Ext_R^1(N,M)=0$ is projective, or equivalently, free.
The residue field $k$ is a typical example of a strong test module for projectivity.
Ramras shows that the maximal ideal $\m$ is a strong test module for projectivity. 
He also proves that every strong test module for projectivity has depth at most one.
Using the rigidity theorem for Tor modules, Jothilingam \cite{J1} proves that when $R$ is a regular local ring, every $R$-module of depth at most one is a strong test module for projectivity.
Our Theorem \ref{main} yields that the converse of this Jothilingam's result also holds true.

\begin{corollary}\label{stp}
The following seven conditions are equivalent.
\begin{enumerate}[\rm (1)]
\item
$R$ is regular.
\item 
Every $R$-module of depth at most one is a strong test module for projectivity.
\item
Every $R$-module of depth zero is a strong test module for projectivity.
\item
Every $R$-module of depth zero and of finite projective dimension is a strong test module for projectivity.
\item
There exists a strong test $R$-module for projectivity of depth zero and of finite projective dimension.
\item
There exists a strong test $R$-module for projectivity of finite projective dimension.
\item
There exist a semidualizing $R$-module $C$ and a strong test $R$-module for projectivity of finite $C$-projective dimension.
\end{enumerate}
\end{corollary}

\section{Proofs of the results}

Let $M$ be an $R$-module.
Let
$$
\cdots \to F_n \overset{\partial _n}{\to} F_{n-1} \to \cdots \to F_1 \overset{\partial _1}{\to} F_0 \to M \to 0
$$
be a minimal free resolution of $M$.
The $n$th {\em syzygy} $\Omega^nM$ of $M$ is defined as the image of the map $\partial_n$.
We define the {\em (Auslander) transpose} of $M$ to be the cokernel of the map $\Hom_R(\partial_1,R):\Hom_R(F_0,R)\to\Hom_R(F_1,R)$, and denote it by $\tr M$.
Note that the $n$th syzygy and the transpose are uniquely determined up to isomorphism.

To prove our theorem, we establish four lemmas.

\begin{lemma}\label{lem0}
For every $R$-module $M$, there is an isomorphism $\Hom_R(\tr M,R)\cong \Omega^2 M$ of $R$-modules.
\end{lemma}

It follows from \cite[Proposition (2.6)(d) and Line 1 in Page 143]{AB} that $\Hom_R(\tr M,R)$ is isomorphic to $\Omega^2 M$ up to free summands.
Since we define the syzygies and transpose of a module by using its minimal free resolution, we can prove that there is an $R$-isomorphism $\Hom_R(\tr M,R)\cong \Omega^2 M$ in the strict sense.
The proof is standard, and we omit it.

In what follows, let $t=\depth R$.

\begin{lemma}\label{lem1}
Let $C$ be a semidualizing $R$-module, and let $M$ be an $R$-module which is locally free on the punctured spectrum of $R$.
Then:
\begin{enumerate}[\rm (1)]
\item
$\Ext_R^i(C\tensor_RM,C) \cong \Ext_R^i(M,R)$ for all $1\le i\le t$.
\item
There is an injection $\Ext_R^{t+1}(C\tensor_RM,C) \to\Ext_R^{t+1}(M,R)$.
\end{enumerate}
\end{lemma}

\begin{proof}
We have isomorphisms $\RHom_R(C\tensor_R^{\bf L}M,C) \cong \RHom_R(M,\RHom_R(C,C)) \cong \RHom_R(M,R)$ in the derived category of $R$.
From this we obtain a spectral sequence
$$
\E_2^{pq}=\Ext_R^p(\Tor^R_q(C,M),C) \Rightarrow \H^{p+q}=\Ext_R^{p+q}(M,R).
$$
By \cite[Page 63(1)]{G} we have $\depth_RC=t$.
Since $\Tor^R_q(C,M)$ has finite length for $q>0$, we have $\E_2^{pq}=0$ when $q>0$ and $p<t$.
Thus we obtain an isomorphism $\H^i\cong\E_2^{i0}$ for $1\le i\le t$, and an injective homomorphism $\E_2^{t+1,0} \to \H^{t+1}$.
\end{proof}

\begin{lemma}\label{prop1}
Let $C$ be a semidualizing $R$-module.
\begin{enumerate}[\rm(1)]
\item
If $\Ext_R^i(X,C)=0$ for $1\le i\le t+1$, then $X$ belongs to $\Xc_C(R)$.
\item
The module $C\tensor_R\tr\Omega^{t+1}k$ belongs to $\Xc_C(R)$.
\end{enumerate}
\end{lemma}

\begin{proof}
(1) Let $M$ be an $R$-module of finite $C$-projective dimension.
Then we have $\cpd_RM = \pd_R(\Hom_R(C,M)) \le \depth R =t$ by \cite[Theorem 2.11(c)]{TW}, and there is an exact sequence
$$
0 \to C^{b_t} \to C^{b_{t-1}} \to \cdots \to C^{b_1} \to C^{b_0} \to M \to 0,
$$
where each $b_i$ is a nonnegative integer.
We have short exact sequences
$$
0 \to M_{i+1} \to C^{b_i} \to M_i \to 0 \quad (0 \leq i \leq t-1),
$$
where $M_0=M$ and $M_t=C^{b_t}$.
Since $\Ext_R^i(X,C)=0$ for $1 \leq i \leq t+1$, there are isomorphisms $\Ext_R^1(X,M) = \Ext_R^{1}(X,M_0) \cong \Ext_R^{2}(X,M_1) \cong \cdots \cong \Ext_R^{t+1}(X,M_t) = \Ext_R^{t+1}(X,C^{b_t}) = 0$.
Hence $X$ is in $\Xc_C(R)$.

(2) Note that $\Ext_R^i(V,R)=0$ for every $R$-module $V$ annihilated by $\m$ and every integer $i\le t-1$.
Hence $\Ext_R^i(\Ext_R^j(k,R),R)=0$ for $1 \leq j \leq t+1$ and $0 \leq i \leq j-2$.
It follows from \cite[Proposition (2.26)]{AB} that the $R$-module $\Omega ^{t+1}k$ is {\em $(t+1)$-torsionfree}, that is, $\Ext_R^i(\tr\Omega^{t+1}k,R)=0$ for all $1\le i\le t+1$.
Since $\tr\Omega^{t+1}k$ is locally free on the punctured spectrum of $R$, Lemma \ref{lem1} implies that for $1\le i\le t+1$ there exists an injection $\Ext_R^i(C\tensor_R\tr\Omega^{t+1}k,C) \to \Ext_R^i(\tr\Omega^{t+1}k,R)$.
Therefore $\Ext_R^i(C\tensor_R\tr\Omega^{t+1}k,C)=0$, and thus $C\tensor_R\tr\Omega^{t+1}k$ is in $\Xc_C(R)$ by (1).
\end{proof}

Now let us verify that the inclusion relations stated in the introduction hold.

\begin{lemma}\label{prop2}
\begin{enumerate}[\rm(1)]
\item
Let $C$ be a semidualizing $R$-module.
Then
$$
\Xc_{C}(R)\supset\Gc_{C}(R)\supset\add C,\ \add R.
$$
\item
One has
$$
\Xc(R)\supset\Gc(R)\supset\Ic(R)\supset\add R.
$$
\end{enumerate}
\end{lemma}

\begin{proof}
(1) It follows from Lemma \ref{prop1}(1) that $\Xc_C(R)$ contains $\Gc_C(R)$.
Since $R$ and $C$ belong to $\Gc_C(R)$, the subcategories $\add R$ and $\add C$ are contained in $\Gc_C(R)$.

(2) Considering in (1) the semidualizing module $R$, we see that $\Xc(R)$ contains $\Gc(R)$.
The other two inclusion relations are obtained by \cite[Theorem (1.4)]{AGP}.
\end{proof}

Now we can prove the theorem.

\begin{tpf}
(1) (a) $\Rightarrow$ (b):
Let $C$ be a dualizing $R$-module.
By \cite[Corollary (3.9)]{S}, $R$ is a Cohen-Macaulay local ring.
Denote by $\CM(R)$ the full subcategory of $\mod R$ consisting of all maximal Cohen-Macaulay $R$-modules.
Note here that $\CM (R)$ coincides with $\Gc_C(R)$.
By virtue of the Cohen-Macaulay approximation theorem \cite[Theorem 1.1]{ABu}, for each $R$-module $M$ there exists an exact sequence
$$
0 \to Y \to X \to M \to 0
$$
of $R$-modules such that $X$ is maximal Cohen-Macaulay and $Y$ has finite $C$-projective dimension.
Now suppose that $M$ belongs to $\Xc_C(R)$.
Then the above exact sequence splits because it corresponds to an element of $\Ext_R^1(M,Y)$, which vanishes as $M\in\Xc_C(R)$ and $\cpd_RY<\infty$.
Hence $M$ is isomorphic to a direct summand of $X$.
In particular, $M$ is a maximal Cohen-Macaulay $R$-module, and thus $\Xc_C(R)$ is contained in $\CM (R)=\Gc_C(R)$.
Lemma \ref{prop2}(1) implies that $\Xc_C(R)$ coincides with $\Gc_C(R)$.

(b) $\Rightarrow$ (a):
Lemma \ref{prop1}(2) shows that $C\tensor_R\tr\Omega^{t+1}k$ belongs to $\Xc_C(R)=\Gc_C(R)$.
By definition, the $C$-dual module $\Hom_R(C\tensor_R\tr\Omega^{t+1}k,C)$ is also in $\Gc_C(R)$.
There are isomorphisms
\begin{align*}
\Hom_R(C\tensor_R\tr\Omega^{t+1}k,C) & \cong \Hom_R(\tr\Omega^{t+1}k,\Hom_R(C,C)) \\
& \cong \Hom_R(\tr\Omega^{t+1}k,R) \\
& \cong \Omega^2 \Omega^{t+1}k =\Omega^{t+3}k
\end{align*}
by Lemma \ref{lem0}.
Hence $\Omega^{t+3}k$ belongs to $\Gc_C(R)$, and $\Ext^{i+t+3}_R(k,C)\cong\Ext^i_R(\Omega^{t+3}k,C)=0$ for all $i>0$, which shows that $C$ has finite injective dimension.
Thus $C$ is a dualizing $R$-module.

(2) Letting $C=R$ in (1), we observe that the assertion holds.

(3) (a) $\Rightarrow$ (b):
Since $R$ is Gorenstein, we have $\Xc(R)=\Gc(R)=\CM (R)$ by (2).
We see from \cite[Theorem (1.4)]{AGP} that $\Ic(R)$ now coincides with $\CM(R)$.

(b) $\Rightarrow$ (a):
The assumption and Lemma \ref{prop2}(2) give the equalities $\Xc(R)=\Gc(R)=\Ic(R)$, and $R$ is Gorenstein by (2).
Hence we have $\Ic(R)=\Gc(R)=\CM(R)$.
As $\Omega^nk$ belongs to $\CM(R)=\Ic(R)$ for $n\gg0$, the local ring $R$ is a complete intersection by \cite[Lemma (1.9) and Theorem (1.3)]{AGP}.

(4) (a) $\Rightarrow$ (b):
If $R$ is regular, then it is Gorenstein and every maximal Cohen-Macaulay $R$-module is free.
By (2) we have $\Xc(R)=\Gc(R)=\CM (R)=\add R$.

(b) $\Rightarrow$ (c) and (b) $\Rightarrow$ (d):
Let $C=R$.

(d) $\Rightarrow$ (b):
The module $C$ belongs to $\Xc_C(R)=\add R$ by Lemma \ref{prop1}(1), so it is isomorphic to $R^r$ for some $r\ge0$.
Since $C$ is semidualizing, we have $R \cong \Hom_R(C,C) \cong R^{r^2}$.
Hence we have $r=1$ and $C \cong R$.
Therefore it holds that $\Xc(R)=\add R$.

(c) $\Rightarrow$ (a):
Lemma \ref{prop1}(2) says that $C \tensor_R \tr \Omega^{t+1}k$ is in $\Xc_C(R)=\add C$.
Similarly to the proof of (b) $\Rightarrow$ (a) in (1), we observe that $\Omega^{t+3}k$ is isomorphic to $\Hom_R(C \tensor_R \tr \Omega^{t+1}k,C)$, which is in $\add R$ because for a module $M\in\add C$ the $C$-dual module $\Hom_R(M,C)$ is in $\add R$.
Thus the $R$-module $k$ has finite projective dimension, and hence $R$ is regular.
\qed
\end{tpf}

Now let us give a proof of the corollary.

\begin{spf}
(1) $\Rightarrow$ (2):
This implication follows from \cite[Theorem 1]{J1}.

(2) $\Rightarrow$ (3) $\Rightarrow$ (4) and (5) $\Rightarrow$ (6):
These implications are obvious.

(4) $\Rightarrow$ (5):
Take a maximal $R$-regular sequence $x_1, x_2, \ldots , x_t$.
Then $R/(x_1, x_2, \ldots , x_t)$ is an $R$-module of depth zero and of finite projective dimension.

(6) $\Rightarrow$ (7):
Letting $C=R$ shows this implication.

(7) $\Rightarrow$ (1):
Let $M$ be a strong test $R$-module for projectivity with $\cpd_RM<\infty$.
Let $N$ be a module in $\Xc_C(R)$.
Then we have $\Ext_R^1(N,M)=0$.
Since $M$ is a strong test module for projectivity, $N$ is a free $R$-module.
Thus $\Xc_C(R)$ is contained in $\add R$.
Lemma \ref{prop2}(1) implies $\Xc_C(R)=\add R$, and $R$ is regular by Theorem \ref{main}(4).
\qed
\end{spf}

\begin{ac}
The authors are grateful to Naoya Hiramatsu and Takeshi Yoshizawa for their valuable comments.
They also thank the referee for his/her careful reading.
\end{ac}

\end{document}